\documentclass[a4paper,10pt]{article}
\usepackage{amsfonts}
\usepackage{amsmath}
\usepackage{amssymb}
\hoffset= - 1.0in         

\begin{document}
\hfill ITEP/TH-40/07

\bigskip

\centerline{ \Large{Two-Body systems from SL(2,$\mathbb{C}$)-tops}}

\bigskip

\centerline{\it A.Smirnov {\footnote{E-mail: asmirnov@itep.ru}}}

\bigskip

\centerline{ITEP, Moscow, Russia}

\bigskip

It is shown that sl(2,$\mathbb{C}$) Euler-Arnold tops are equivalent to the two-body
systems of Calogero-Moser type. We prove  that generic Hamiltonians of sl(2,$\mathbb{C}$) tops are
equivalent to one of three canonical Hamiltonians. For all canonical
Hamiltonians the corresponding  two-body system is found. Bosonisation formulas for
each case are obtained explicitly. Relations  with Antonov-Zabrodin-Hasegawa
R-matrix are discussed.
\section{Introduction}
A most famous examples of sl(2,$\mathbb{C}$) Euler-Arnold tops are the so called $XYZ$, $XXZ$, $XXX$ tops.
They can be defined by their Lax operators.
The Hamiltonian of the system is identified with free term in the expansion of $Tr(L(z)^2)$ into series
with respect to spectral parameter $z$.
The Lax operators $L(z)$ for these systems can be found from the equation for the linear Poisson brackets:
$$ \{L(z)\otimes1,1\otimes L(z)\}=[r,1\otimes L(z)+L(z)\otimes 1],$$
where r is the elliptical, the trigonometrical or the rational classical Baxter sl(2,$\mathbb{C}$) r-matrix.\\
\indent
In the framework of the Hitchin approach to integrable systems \cite{Hit}, \cite{Gorsky},
\cite{GorskyNekrasov}, \cite{Nekrasov} the connection between the $XYZ$ top,
and the two-body Calogero-Moser system was established in \cite{LOZ}. It was shown that the Lax operators
$L_{ell}^{top}(z)$ and $L^{CM}_{ell}(z)$  for systems are conjugated by the gauge transformation:
\begin{equation}
\label{SHC}
L^{top}_{ell}(z)=\Xi(z)L^{CM}_{ell}(z)\Xi^{-1}(z).
\end{equation}
Here the  matrix $\Xi(z)$ is singular at the point $z=0$.
Let
$$
M^{top}=\lbrace S_{j}\ \ \ j=1,2,3\,:S_{1}^2+S_{2}^2+S_{3}^2=\nu^2 \rbrace
$$
be a sl(2,$\mathbb{C}$) coadjoint orbit. It is a phase space for the sl(2,$\mathbb{C}$) Euler-Arnold tops
with a Poisson bracket $\{S_{i},S_{j}\}=2\,i\,\epsilon_{ijk}\,S_{k}$. Define a coordinate $q$
and a momentum $p$ in the center mass of two-body Calogero-Moser system. The so called bosonisation
formulas can be read off from (\ref{SHC}). These formulas explicitly give the symplectomorphism between systems:
\begin{eqnarray}
\label{bf1}
&& S_{1}=-\frac{\theta_{10}(0)\theta_{10}(2q)}{\vartheta^{'}(0)\vartheta(2q)}p+
 \frac{ \theta_{10}^2(0)\theta_{00}(2q)\theta_{01}(2q)}{\theta_{00}(0)\theta_{01}(0) \vartheta^2(2q)}\nu,\\
\label{bf2}
&&S_{2}=\frac{\theta_{00}(0)\theta_{00}(2q)}{i\vartheta^{'}(0)\vartheta(2q)} p-
\frac{ \theta_{00}^2(0)\theta_{10}(2q)\theta_{01}(2q)}{i\theta_{10}(0)\theta_{01}(0 )\vartheta^2(2q)}\nu,\\
\label{bf3}
&&S_{3}=-\frac{\theta_{01}(0)\theta_{01}(2q)}{\vartheta^{'}(0)\vartheta(2q)} p+
\frac {\theta_{01}^2(0)\theta_{00}(2q)\theta_{10}(2q)}{\theta_{00}(0)\theta_{10}(0 )\vartheta^2(2q)}\nu,
\end{eqnarray}
where $\vartheta(q)$, $ \theta_{\alpha \beta}(q)  $ are the standard elliptic theta function and the
theta functions with characteristics respectively.
Evidently, the  connection between the tops and the Calogero-Moser
systems should exist in trigonometrical and rational limits.
However, in the naive trigonometric limit $\tau\rightarrow
i\infty$ the matrix elements of $\Xi(z)$ are divergent. This difficulty was overcame in \cite{AHZ}
by applying an additional gauge transformation $A(\tau)$ depending on modular  parameter $\tau$ which
 is singular at point $\tau=i\infty$:
$$
A(\tau)\,L^{top}_{ell}(z)\,A^{-1}(\tau)=A(\tau)\,\Xi\, L^{CM}_{ell}(z)\,\Xi^{-1}\,A^{-1}(\tau)
$$
Upon taking the limits:
$$
L_{trig}^{' top}(z)=\lim_{\tau\rightarrow i
\infty}A(\tau)\,L_{ell}^{top}(z)\,A^{-1}(\tau),\ \ \
\Xi^{'}(z)=\lim_{\tau\rightarrow i \infty}A(\tau)\,\Xi(z),\ \ \
 L^{CM}_{trig}(z)=\lim_{\tau\rightarrow i \infty}\,L^{CM}_{ell}(z),
$$
we get the regular trigonometric version of (\ref{SHC}):
\begin{equation} \\
\label{NewLimit}
L^{' top}_{trig}(z)=\Xi^{'}(z)L^{CM}_{trig}(z)\Xi^{' -1}(z).
\end{equation}
Due to the fact that matrix $A(\tau)$ is singular at the point $\tau=i
\infty$ this new trigonometric Lax operator is not gauge equivalent to the standard one:
$$
L^{top}_{trig}(z)=\lim_{\tau\rightarrow i\infty}L_{ell}^{top}(z),
$$
In this way we obtain a distinct from the ordinary $XXZ$-top Hamiltonian system on the coadjoint sl(2,$\mathbb{C}$) orbit.
 Therefore,
in addition to standard $XXZ$ and $XXX$ tops there exist nonequivalent to them trigonometrical
and rational degenerations of elliptic $XYZ$ top. It follows from (\ref{NewLimit}) that these tops describe the
 evolution of the
trigonometrical and the rational Calogero-Moser systems. In what follows we will denote these new tops
as $XXZ^{'}$ and $XXX^{'}$ tops respectively.\\
\indent
Our approach to the problem is different from \cite{AHZ}. In section 3 we show that a general
 Hamiltonian of sl(2,$\mathbb{C}$) top is gauge
equivalent to one of three canonical Hamiltonian presented below. We obtain the bosonisation formulas
 for all top Hamiltonians in a simple way in section 4.
By these means in sections 5-7 we receive the bosonisation formulas for each canonical
Hamiltonian. We see that the bosonisations corresponding to the canonical Hamiltonians
map tops to the two-body Calogero-Moser systems. Therefore presented canonical
Hamiltonians are the Hamiltonians of  $XYZ$, $XXZ^{'}$ and $XXX^{'}$ tops.
The bosonisation formulas obtained for $XYZ$ top coincide with ones
obtained in \cite{LOZ}. Formulas received for $XXZ^{'}$ and $XXX^{'}$ tops are
their natural degenerations.\\
\indent
In section 8 we show that in process of obtaining trigonometric bosonisation formulas
from elliptical ones (\ref{bf1})-(\ref{bf3}) by tending $\tau\rightarrow i\infty$
we obtain divergent results. To get finite expressions
we perform regularization on the algebra sl(2,$\mathbb{C}$) by introducing
singular gauge transformation depending on $\tau$. We see that this transformation is the
classical version of singular gauge transformation first introduced in \cite{GorskyZabrodin},
\cite{AHZ}.
\section{Tops on Lie algebras}
Let $\mathfrak{G}$ be a semisimple Lie algebra. The Euler-Arnold top on the algebra $\mathfrak{G}$ is the
 Hamiltonian system on the coadjoint orbit:
\begin{equation}
\label{eq1}
M^{rot}=\lbrace S\in \mathfrak{G}^{*}:S=g^{-1}S_{0}g \rbrace,
\end{equation}
where $g$ are elements of corresponding group $g\in G$.
The phase space is equipped with a non-degenerate Kirillov-Kostant symplectic form:
\begin{equation}
\label{eq2}
w=<S_{0},dgg^{-1}\wedge dgg^{-1}>.
\end{equation}
Where $<,>$ is the Killing form. Dynamics of the system is described by the Hamiltonian:
\begin{equation}
\label{eq3} \\
H_{J}=J_{ij}S_{i}S_{j},
\end{equation}
where $J$ is some quadratic form on the algebra $\mathfrak{G}$. It is defined up to the quadratic Casimir element:
\begin{equation}
\Omega=\Omega_{ij}S_{i}S_{j},\ \ \ \lbrace \Omega,A\rbrace=0\ \ \forall A\in \mathfrak{G}^{*}.
\end{equation}
The Poisson structure on the algebra of functions can be described explicitly. On linear functions
the Poisson brackets coincide with commutator:
$$
\lbrace A , B \rbrace = [A,B], \ \ \ \  A,B \in \mathfrak{G}^{*},
$$
whereupon Poisson brackets on polynomial functions can be evaluated by the Leibniz rule:
$$
\lbrace AB,C\rbrace=\lbrace A , C\rbrace B + \lbrace B , C\rbrace A,\ \ \ \ \ \ \ \ \ \
\ \ \ \ \ \ \ \ \ A,B,C\in F(M^{rot}).
$$
Hence the equations of motion of the top in some basis $S_{i}$:
\begin{equation}
\label{eq4} \\
\frac{dS_{k}}{dt}=\lbrace H_{J},S_{k} \rbrace = \lbrace J_{ij} S_{i} S_{j},S_{k} \rbrace=2J_{ij}C_{mjk}S_{i}S_{m},
\end{equation}
where $ C_{i j k} $ is the tensor of structure constants in this basis.\\
\indent
Hereby, the top is defined by choosing the quadratic form on the algebra. However two different
quadratic form $J_{1}$ and $J_{2}$ can lead to the same equation of motion (\ref{eq4}).
If it is the case we will say that forms $J_{1}$ and $J_{2}$ are equivalent.
More precisely\\
\textbf{Proposition.}  \textit{Let two quadratic forms are equivalent with respect to
the action of a the automorphism group of the algebra $\mathfrak{G}$, i.e. $J_{2ij}=J_{1km}T_{ki}T_{mj}$,
where $T_{ij}$ $\in$ Aut($\mathfrak{G}$), then they define the same equations of motions.}\\
\textit{Prove.} For $H_{J_{1}}$ from (\ref{eq4}) we have:
\begin{equation}
\label{eq5}
\frac{dS_{k}}{dt}=2J_{1ij}C_{mjk}S_{i}S_{m},
\end{equation}
for $H_{J_{2}}$:
$$
\frac{dS_{k}^{'}}{dt}=\lbrace H_{J_{2}},S_{k}^{'}\rbrace=\lbrace J_{2lm}S_{l}S_{m},S_{k}^{'} \rbrace=
\lbrace J_{1ij}T_{il}T_{jm}S_{l}S_{m},S_{k}^{'}\rbrace=\lbrace J_{1ij}S_{i}^{'}S_{j}^{'},S_{k}^{'}\rbrace,
$$
where we denote $S_{i}^{'}=T_{ij}S^{j}$. Therefore:
\begin{equation}
\label{eq6}
\frac{dS_{k}^{'}}{dt}=2J_{1ij}C_{jkm}^{'}S_{i}^{'}S_{m}^{'},
\end{equation}
where $ C_{ijk}^{'}$ is the tensor of structure constants in the basis $S_{i}^{'}$. The action of the
automorphism group preserves the commutation relations $ \lbrace S_{i}^{'},S_{j}^{'} \rbrace=C_{ijk}S_{k}^{'}$, and
 hence $C_{ijk}^{'}=C_{ijk}$. Hereby the equations (\ref{eq5}) and (\ref{eq6}) coincide.\\
It follows from above that the space of the quadratic forms on the algebra $\mathfrak{G}$ fall into equivalence classes.
The equivalent quadratic forms correspond to the equivalent tops.
\section{sl(2,$\mathbb{C}$)-case}
It follows from section 2 that the classification of the tops on the algebra $\mathfrak{G}$ is reduced
the to the classification of equivalence classes of quadratic forms on the algebra with respect to the action of the
 automorphism group. It turn out that in the case $\mathfrak{G}$=sl(2,$\mathbb{C}$) the matrix of any
quadratic form can be represented in one of the three canonical form. The canonical forms of
matrix in some basis with the relations:
\begin{equation}
\label{stbas}
\lbrace S_{i},S_{j} \rbrace=2i\epsilon_{ijk}S_{k},
\end{equation}
are given by matrices:
\begin{equation}
\label{eq7}
XYZ)\\\
J=\left[\begin{array}{ccc}
\alpha&0&0\\0&\beta&0\\0&0&\gamma
\end{array}\right],\\\
XXZ^{'} )\\\
J=\left[\begin{array}{ccc}
\alpha&i\alpha&0\\i\alpha&-\alpha&0\\0&0&\beta
\end{array}\right],
\\\
XXX^{'})\\\
J=\left[\begin{array}{ccc}
\alpha&i\alpha&\beta\\i\alpha&-\alpha&i\beta\\\beta&i\beta&0
\end{array}\right].\\\
\end{equation}
We will call bases of this kind the standard bases. 
In the standard basis (\ref{stbas}) the Casimir function take the form:
\begin{equation} \\
\label{cas}
\Omega=S_{1}^2+S_{2}^2+S_{3}^2.
\end{equation}
The Casimir (\ref{cas}) is invariant with respect to the action of automorphism group. Thus,
the matrix of $A\,\in\,Aut(sl(2,\mathbb{C}))$ is orthogonal in the standard basis $A A^{t}=1$. Therefore,
the the matrix of quadratic form is tramsformed as the matrix of operator:
$$
J^{'}=A J A^{t}=A J A^{-1}.
$$
So, the the canonical forms of the matrices of a quadratic forms on the algebra is exhausted by the symmetric matrices with different Jordan structures~(\ref{eq7}).

The Hamiltonians corresponding to these matrices have the forms:
\begin{eqnarray}
&&XYZ)\,H=\alpha S_{1}^2+\beta S_{2}^2+\gamma S_{3}^{2},\\
&&XXZ^{'})\,H=\alpha (S_{1}+i S_{2})^2+\beta S_{3}^2,\\
&&XXX^{'})\,H=\alpha (S_{1}+i S_{2})^2+\beta S_{3}(S_{1}+i S_{2});
\end{eqnarray}
It is useful to represent these Hamiltonians in the Chevalley basis:
\begin{eqnarray}
\label{c1}
&& XYZ)\,H=\alpha (e+if)^2+\beta (e-if)^2+\gamma h^2 \\
\label{c2}
&& XXZ^{'})\,H=\alpha e^2+\beta h^2\\
\label{c3}
&& XXX^{'})\,H=\alpha e^2+\beta eh
\end{eqnarray}
\begin{equation}
\label{eq9}
\lbrace h , e \rbrace = 2e,\\\
\lbrace h , f \rbrace =-2f,\\\
\lbrace e , f \rbrace =h.
\end{equation}
Therefore general sl(2,$\mathbb{C}$)-top is described by one of the three Hamiltonians
(\ref{c1})-(\ref{c3}). In next sections we will investigate these Hamiltonians
and find the corresponding two-particle systems. Thus, we present the complete connection between
the sl(2,$\mathbb{C}$)-tops and the systems of two particles. 
\section{Particles-Tops correspondence}
In this section we describe the main technique for obtaining the bosonisation formulas. The system of two particles is
 described by the Hamiltonian:
\begin{equation} \\
\label{eq10}
H=p^2+U(q).
\end{equation}
It is defined on the two-dimensional phase space :
\begin{equation}
\label{eq11}
M^{part}=\lbrace (p,q) \in \mathbb{C}^2\rbrace , \\\
w=dp\wedge dq.
\end{equation}
with the canonical Poisson bracket:
\begin{equation} \\
\label{eq12}
\lbrace f,g \rbrace = \frac{\partial f}{\partial p}\frac{\partial g}{\partial q}-
\frac{\partial f}{\partial q}\frac{\partial g}{\partial p} , \\\
\lbrace p,q \rbrace=1 . \\
\end{equation}\\
\indent
Consider a top with quadratic form $J$. To define the corresponding system of two particles we need  explicit
expressions for generators of the algebra $S_{i}(p,q)$ in terms of $p$ and $q$. With respect to the Poisson bracket
 (\ref{eq12}) these functions obey
the standard commutation relations (\ref{stbas}):
$$
\lbrace S_{i}(p,q),S_{j}(p,q) \rbrace=2i\epsilon_{ijk}S_{k}(p,q)
$$
Assume that there exist a map $(p,q)\rightarrow S_{k}$ that transform the Hamiltonian (\ref{eq3}) into (\ref{eq10}).
 Then it is necessary
to consider a linear in $p$ anzats for  the generators of the algebra (the bosonisation formulas):
\begin{equation}
\label{eq13}
h=f_{h}(q)p+g_{h}(q),\ \
e=f_{e}(q)p+g_{e}(q),\ \
f=f_{f}(q)p+g_{f}(q),
\end{equation}
where function $f_{h}(q),f_{e}(q),f_{f}(q),g_{h}(q),g_{e}(q),g_{f}(q)$ should be defined. Commutation
relations (\ref{eq9}) lead to the system of differential equations:
\begin{eqnarray}
\label{st1}
&&f_{h}(q)\dfrac{df_{f}(q)}{dq}-f_{f}(q)\dfrac{df_{h}(q)}{dq}=-2f_{f}(q),\\
&&f_{h}(q)\dfrac{dg_{f}(q)}{dq}-f_{f}(q)\dfrac{dg_{h}(q)}{dq}=-2g_{f}(q),\\
&&f_{h}(q)\dfrac{df_{e}(q)}{dq}-f_{e}(q)\dfrac{df_{h}(q)}{dq}=2f_{e}(q),\\
&&f_{h}(q)\dfrac{dg_{e}(q)}{dq}-f_{e}(q)\dfrac{dg_{h}(q)}{dq}=2g_{e}(q),\\
&&f_{e}(q)\dfrac{df_{f}(q)}{dq}-f_{f}(q)\dfrac{df_{e}(q)}{dq}=f_{h}(q),\\
\label{st2}
&&f_{e}(q)\dfrac{dg_{f}(q)}{dq}-f_{f}(q)\dfrac{dg_{e}(q)}{dq}=g_{h}(q).
\end{eqnarray}
Using (\ref{eq13}) we rewrite the Hamiltonian of the top (\ref{eq3}). It has the form of quadratic in p
 polynom:
\begin{equation}
\label{eq15}
H=\Lambda_{1}(f_{h} , f_{e}, f_{f})p^2+\Lambda_{2}(f_{h},f_{e},f_{f},g_{h},g_{e},g_{f})p+\Lambda_{3}(g_{h },g_{e},g_{f}),
\end{equation}
where particular forms of the coefficients are defined by $J$. By comparing (\ref{eq15}) with (\ref{eq10})
we find:
\begin{eqnarray}
\label{bc1}
&&\Lambda_{1}(f_{h} , f_{e} ,f_{f})=1,\\
\label{bc2}
&&\Lambda_{2}(f_{h},f_{e},f_{f},g_{h},g_{e},g_{f})=0.
\end{eqnarray}
There is an additional condition, because the Casimir element is a constant on orbits:
\begin{equation}
\label{eq17}
h^2+4ef=\nu^2.
\end{equation}
The last equations (\ref{bc1})-(\ref{eq17}) play the role of boundary conditions for the system  (\ref{st1})-(\ref{st2})
 and  allow us to solve this system explicitly. As a result we  obtain the function $\Lambda_{3}(q)$ that play the role of
the interaction potential.
\section{The $XXX^{'}$-case}
Let us consider the system of two particles corresponding to $XXX^{'}$ top following the receipt of the
previous section. According with (\ref{c3}) the Hamiltonian of top in $XXX^{'}$-case has the following form:
$$
H=\alpha e^2+\beta eh.
$$
where we assume $\alpha=1$ for normalization:
\begin{equation}
\label{eq18}
H=e^2+\beta eh.
\end{equation}
Substituting (\ref{eq13}) into (\ref{eq18}) we find the following boundary conditions:
\begin{eqnarray}
\label{xxxbc1}
&&\Lambda_{1}(f_{e},f_{h},f_{f})=f_{e}^{2}(q)+\beta f_{e}(q)f_{h}(q)=1,\\
\label{xxxbc2}
&&\Lambda_{2}(f_{e},f_{h},f_{f},g_{e},g_{h},g_{f})=2 f_{e}(q)g_{e}(q)+\beta f_{e}(q)g_{h}(q)+\beta f_{h}(q)g_{e}(q)=0.
\end{eqnarray}
\textbf{Proposition.} \textit{The solution of the system (\ref{st1})-(\ref{st2}) satisfying
the boundary conditions (\ref{xxxbc1}),(\ref{xxxbc2}), and (\ref{eq17}) has the following form:}
\begin{eqnarray}
\label{bosxxx}
&&f_{h}(q)=\frac{1-q^2\beta^2}{q\beta^2},\\
&&g_{h}(q)=\frac{\nu(q^2\beta^2+1)}{2q^2\beta^2}.\\
&&f_{e}(q)=-\frac{1}{q\beta},\\
&&g_{e}(q)=-\frac{\nu}{2q^2\beta},\\
&&f_{f}(q)=\frac{(1-q^2\beta^2)^2}{4q\beta^3},\\
&&g_{f}(q)=\frac{\nu(1-q^2\beta^2)(1+3q^2\beta^2)}{8q^2\beta^3}.
\end{eqnarray}
Hereby the bosonisation formulas have the form:
\begin{eqnarray}
&&e=-\frac{p}{q\beta}-\frac{\nu}{2q^2\beta},\\
&&h=\frac{1-q^2\beta^2}{q\beta^2}p+\frac{\nu(q^2\beta^2+1)}{2q^2\beta^2},\\
&&f=\frac{(1-q^2\beta^2)^2}{4q\beta^3}p+\frac{\nu(1-q^2\beta^2)(1+3q^2\beta^ 2) }{8q^2\beta^3}.
\end{eqnarray}
For the Hamiltonian we have:
\begin{equation}
\label{eq23}
H=e^2+\beta eh=p^2-\frac{\nu^2}{(2q)^2}.
\end{equation}
This Hamiltonian describe the two-particle Calogero-Moser system with a rational potential.
 Thus, we have proved that $XXX^{'}$ top is  equivalent to the rational Calogero-Moser system.
\section{The  $XXZ^{'}$-case}
Let us repeat previous calculations for $XXZ^{'}$ top. According with (\ref{c2}) in this case the Hamiltonian of top
 has the form:
$$
H=e^2+\gamma^2 h^2,
$$
where we assume $\alpha=1$, $\beta=\gamma^2$ for simplicity. The system of the boundary
conditions (\ref{bc1}), (\ref{bc2}) in this case have the form:
\begin{eqnarray}
\label{r1}
&&\Lambda_{1}(f_{e},f_{h},f_{f})=f_{e}^2(q)+\gamma^2 f_{h}^2(q)=1,\\
\label{r2}
&&\Lambda_{2}(f_{e},f_{h},f_{f},g_{e},g_{h},g_{f})=f_{e}(q)g_{e}(q)+\gamma^2 f_{h}(q)g_{h}(q)=0.
\end{eqnarray}
\textbf{Proposition.} \textit{Solution of the system (\ref{st1})-(\ref{st2}) satisfying the
boundary conditions (\ref{r1}),(\ref{r2}), and (\ref{eq17}) have the following form:}
\begin{eqnarray}
&&f_{h}(q)=-\dfrac{1}{\gamma th(2\gamma q)},\\
&&g_{h}(q)=-\dfrac{\nu}{sh^2(2\gamma q)},\\
&&f_{e}(q)=-\dfrac{i}{sh(2\gamma q)},\\
&&g_{e}(q)=-\dfrac{i\nu \gamma ch(2\gamma q)}{sh^2(2\gamma q)},\\
&&f_{f}(q)=-\dfrac{i ch^2(2\gamma q)}{4\gamma^2 sh(2\gamma q)},\\
&&g_{f}(q)=\dfrac{i\nu ch(2 \gamma q)(ch^2(2\gamma q)-2)}{4\gamma sh^2(2\gamma q)}.
\end{eqnarray}
The bosonisation formulas (\ref{eq13}) in $XXZ^{'}$-case take the form:
\begin{eqnarray}
\label{bftr1}
&&e=-\frac{i}{sh(2\gamma q)}p-\frac{i\nu \gamma ch(2\gamma q)}{sh^2(2\gamma q)},\\
&&h=-\frac{1}{\gamma th(2\gamma q)}p-\frac{\nu}{sh^2(2\gamma q)},\\
\label{bftr3} &&f=-\frac{i ch^2(2\gamma q)}{4\gamma^2 sh(2\gamma
q)}p+\frac{i\nu ch(2 \gamma q)(ch^2(2\gamma q)-2)}{4\gamma
sh^2(2\gamma q)}.
\end{eqnarray}
 We obtain the following expression for the Hamiltonian:
\begin{equation}
\label{eq28}
H=e^2+\gamma^2h^2=p^2-\frac{\gamma^2\nu^2}{sh^2(2\gamma q)}.
\end{equation}
It corresponds to two-particle Calogero-Moser system with the trigonometric potential.
Thus, we have proved the equivalence of the $XXZ^{'}$-top and the  trigonometric Calogero-Moser system
for two particles.
\section{The $XYZ$-case}
According with (\ref{c1}) the Hamiltonian of $XYZ$ top is given by the expression:
$$
H=\alpha(e+f)^2+\beta(e-f)^2+\gamma h^2.
$$
Let us consider the Casimir function:
$$
\Omega=h^2+4ef=h^2+(e+f)^2-(e-f)^2,
$$
If we subtract $\gamma \Omega $ from the Hamiltonian, and assume $\beta+\gamma=k^2$,$\gamma-\alpha=1$
for normalization we find:
\begin{equation}
\label{eq29}
H=k^2(e-f)^2-(e+f)^2.
\end{equation}
The boundary conditions (\ref{st1})-(\ref{st2}) for this Hamiltonian take the form:
\begin{eqnarray}
\label{r4}
&&\Lambda_{2}=k^2(f_{e}(q)-f_{f}(q))(g_{e}(q)-g_{f} (q))-(f_{e}(q)+f_{f}(q))(g_{e}(q)+g_{f}(q))=0,\\
\label{r5}
&&\Lambda_{1}=k^2(f_{e}(q)-f_{f}(q))^2-(f_{e}(q)+f_{f}(q))^2=1.
\end{eqnarray}
\textbf{Proposition.} \textit{The solution of the system (\ref{st1})-(\ref{st2}) satisfying the
boundary conditions (\ref{r4}),(\ref{r5}), and (\ref{eq17}) has the following form:}
\begin{eqnarray}
&&f_{h}=\dfrac {1}{{sn} \left( 2\,q,k \right) k},\\
&&g_{h}=-{\dfrac {\nu\,{\it dn} \left( 2\,q,k \right) {\it cn}
\left( 2\,q,k \right) }{{\it sn} \left( 2\,q,k \right)
^{2}k}},\\
&&f_{f}=-1/2\,{\dfrac {k{\it cn} \left( 2\,q,k \right) +{\it dn}
\left( 2\,q,k \right) }{{\it sn} \left( 2\,q,k \right) k\sqrt {
1-{k}^{2}}}},\\
&&g_{f}=1/2\,{\dfrac {\nu\, \left( k{\it dn} \left( 2\,q,k \right) +{\it
cn} \left( 2\,q,k \right)  \right) }{ {\it sn}
\left( 2\,q,k \right)  ^{2}k\sqrt {1-{k}^{2}}}},\\
&&f_{e}=-1/2\,{\dfrac {k{\it cn} \left( 2\,q,k \right) -{\it dn}
\left( 2\,q,k \right) }{{\it sn} \left( 2\,q,k \right) k\sqrt {
1-{k}^{2}}}},\\
&&g_{e}=1/2\,{\dfrac {k\nu\,{\it dn} \left( 2\,q,k \right) -\nu\,{\it
cn} \left( 2\,q,k \right) }{  {\it sn} \left( 2\,q,k
\right)  ^{2}k\sqrt {1-{k}^{2}}}},
\end{eqnarray}
where $dn(q,k)$, $sn(q,k)$, $cn(q,k)$ are the  elliptic Jacobi functions with the modular parameter k.
The bosonisation formulas have the form:
\begin{eqnarray}
&&e=-1/2\,{\frac { \left( k{\it cn} \left( 2\,q,k \right) -{\it
dn} \left( 2\,q,k \right)  \right) p}{{\it sn} \left( 2\,q
,k \right) k\sqrt {1-{k}^{2}}}}+1/2\,{\frac {k\nu\,{\it dn}
\left( 2\,q,k \right) -\nu\,{\it cn} \left( 2\,q,k \right) }{
{\it sn} \left( 2\,q,k \right) ^{2}k\sqrt {1-{k
}^{2}}}},\\
&&f=-1/2\,{\frac { \left( k{\it cn} \left( 2\,q,k \right) +{\it
dn} \left( 2\,q,k \right)  \right) p}{{\it sn} \left( 2\,q
,k \right) k\sqrt {1-{k}^{2}}}}+1/2\,{\frac {\nu\, \left( k{\it
dn} \left( 2\,q,k \right) +{\it cn} \left( 2\,q,k \right)
\right) }{ {\it sn} \left( 2\,q,k \right) ^{2}k
\sqrt {1-{k}^{2}}}},\\
&&h={\frac {p}{{\it sn} \left( 2\,q,k \right) k}}-{\frac {\nu\,{\it
dn} \left( 2\,q,k \right) {\it cn} \left( 2\,q,k \right) }
{ {\it sn} \left( 2\,q,k \right) ^{2}k}}.
\end{eqnarray}
It is more useful to present the bosonisation formulas in the standard basis $S_{1}$, $S_{2}$, $S_{3}$:
\begin{equation}
\label{eq33}
S_{1}=e+f,\ \ \ S_{2}=i(f-e),\ \ \ S_{3}=h,\ \ \lbrace S_{i} , S_{j} \rbrace=2i\epsilon_{ijk}S_{k}.
\end{equation}
By expressing the Jacobi functions  in terms of the  elliptic theta functions we obtain the following result:
\begin{eqnarray}
\label{bfell1} && S_{1}={\frac {1}{\sqrt {1-{k}^{2}}}}\left(
-{\frac {{\it cn} \left( 2\,q,k \right) }{{\it sn } \left( 2\,q,k
\right) }}p+{\frac {\nu\,{\it dn} \left( 2\,q,k \right) }{ {\it
sn} \left( 2\,q,k \right)  ^{2}} } \right) =
-\frac{\theta_{10}(0)\theta_{10}(2q)}{\vartheta^{'}(0)\vartheta(2q)}p+\frac{
 \theta_{10}^2(0)\theta_{00}(2q)\theta_{01}(2q)}{\theta_{00}(0)\theta_{01}(0) \vartheta^2(2q)}\nu,\\
\label{bfell2}
&&S_{2}={\frac {1}{i k \sqrt {1-{k}^{2}}}} \left( {\frac {{\it dn} \left( 2\,q,k \right) }{{\it
sn} \left( 2\,q,k \right) }}p-{\frac {\nu\,{\it cn} \left(
2\,q,k \right) }{  {\it sn} \left( 2\,q,k \right)
^{2}}} \right)=
\frac{\theta_{00}(0)\theta_{00}(2q)}{i\vartheta^{'}(0)\vartheta(2q)}p-\frac{
 \theta_{00}^2(0)\theta_{10}(2q)\theta_{01}(2q)}{i\theta_{10}(0)\theta_{01}(0 )\vartheta^2(2q)}\nu,\\
\label{bfell3}
&&S_{3}=\frac{1}{k}\left( {\frac {p}{{\it sn} \left( 2\,q,k \right) }}-{\frac {\nu
\,{\it dn} \left( 2\,q,k \right) {\it cn} \left( 2\,q,k
\right) }{ {\it sn} \left( 2\,q,k \right) ^{2}}
} \right)
=-\frac{\theta_{01}(0)\theta_{01}(2q)}{\vartheta^{'}(0)\vartheta(2q)}p+\frac {
\theta_{01}^2(0)\theta_{00}(2q)\theta_{10}(2q)}{\theta_{00}(0)\theta_{10}(0 )\vartheta^2(2q)}\nu,
\end{eqnarray}
where $\vartheta(q)$, $ \theta_{\alpha \beta}(q)  $ are the standard theta function and theta
 function with characteristics respectively. For the Hamiltonian function we have:
\begin{equation}
\label{eq35}
H=k^2(e-f)^2-(e+f)^2={p}^{2}-\nu^2\wp(2q,k),
\end{equation}
where $\wp(q,k)$is the  Weierstrass elliptic function. The presented Hamiltonian is
the Hamiltonian of the elliptic Calogero-Moser system. It is worth noticing that the presented formulas
coincide with the formulas obtained in \cite{LOZ}.


\section{Connection between XYZ and XXZ$^{'}$ tops }
In this section we explain the interrelation between the  bosonisation formulas corresponding to
$XYZ$ and $XXZ^{'}$ tops. In the limit $k\rightarrow0$ the elliptic Weierstrass function
 behave as:
$$
\wp(q,k)\rightarrow\dfrac{1}{\sin^2(q)},
$$
and the Hamiltonian of elliptic the Calogero-Moser system is  transformed  into the trigonometric one.
Thus, it is natural to suppose that  in this limit $XYZ$ top transform into $XXZ^{'}$ top.
However, the matrices of quadratic forms defining these tops
are not transformed into each other in this way. To understand what is really going on
let us consider the behavior of the bosonisation formulas in the limit $k\rightarrow0$.
From equations (\ref{bfell1})-(\ref{bfell3}) we can find:
\begin{eqnarray}
&&S_{1}\rightarrow-p\dfrac{\cos(2q)}{\sin(2q)}+\nu\frac{1}{\sin^ 2(2q)}\\
&&S_2\rightarrow\frac{1}{i k}\Phi(q,p),\\
&&S_3\rightarrow\frac{1}{k}\Phi(q,p),
\end{eqnarray}
where
\begin{equation}
\Phi(q,p)=p\frac{1}{\sin(2q)}-\nu\frac{\cos(2q)}{\sin^2(2q)}.
\end{equation}
Thus, in this limit the bosonisation formulas diverge. To obtain a finite expressions
we need some regularizing procedure. Let us consider the behavior of Chevalley generators of the
algebra in this limit:
\begin{equation}
h=S_{1},\,e=\dfrac{1}{2}(S_{3}-i\,S_{2})\sim k,\,f=\dfrac{1}{2}(S_{3}+i\,S_{2})\sim\dfrac{1}{k}
\end{equation}
Let us notice that the transformation $h^{'}\rightarrow h$,\, $e^{'}\rightarrow \dfrac{1}{k}e$,\, $f^{'}\rightarrow kf$
 does not change the commutation relations, and we have well defined limit $k\rightarrow 0 $.
After turning back to the standard basis we find the following expressions:
\begin{eqnarray}
\label{trReg1}
&&S_{1}^{'}=h^{'}=- p\frac{\cos(2q)}{4\sin(2q)}+\nu \frac{1}{\sin^{2}
(2q)},\\
\label{trReg2}
&&S_2^{'} =\dfrac{1}{i}(f^{'}-e^{'})= \frac{1}{4 i} \left(p\frac{4 +\cos(2u)}{ \sin(2q)}+
\frac{\nu\,\cos(2q)(5+\sin^{2}(2q))}{\sin^{2}{2q}}\right),\\
\label{trReg3}
&&S_{3}^{'}=(e^{'}+f^{'})=\frac{1}{4}\left(p\frac{4-\cos^2(2q)}{\sin(2q)}-
\frac{\nu\,\cos(2q)(3-\sin^{2}(2q))}{\sin^{2}(2q)}\right).
\end{eqnarray}
It is easy to check that these expressions coincide with (\ref{bftr1})-(\ref{bftr3}) at
some value of constant $\gamma$.  They have correct commutation relations:
\begin{eqnarray}
&&\lbrace S_{i}^{'},S_{j}^{'}\rbrace=2i\epsilon_{ijk}S_{k}^{'},\\
&&(S_{1}^{'})^{2}+ (S_{2}^{'})^{2}+ (S_{3}^{'})^{2}=\nu^{2}.
\end{eqnarray}
In the new basis the matrix of the quadratic form is given by $J^{'}=TJT^{t}$, where $T=A^{-1}BA$ and
\begin{equation}
A=\left[ \begin {array}{ccc} 1&0&0\\\noalign{\medskip}0&1/2\,i&1/2
\\\noalign{\medskip}0&1/2\,i&-1/2\end {array} \right],\,
B=\left[ \begin {array}{ccc} 1&0&0\\\noalign{\medskip}0&{k}^{-1}&0\\ \noalign{\medskip}0&0&k\end {array} \right]
\end{equation}
We see that matrix $T$ is singular at the point $k=0$, but its determinant is equal to one
for any value of parameter $k$. On calculating the matrix $J^{'}$ in the new basis we see that its form
correspond to $XXZ^{'}$ case.
Hereby  $XYZ$ and  $XXZ^{'}$ tops are connected by some singular gauge
transformation,
which depends on modular parameter $k$. Transformations of that kind were introduced
 in \cite{GorskyZabrodin},\cite{AHZ}. In this articles  new solutions
of Yang-Baxter equation for algebras of $A_{n}$ type were obtained. These solutions present non-dynamical R-matrices
for quantum trigonometric Calogero-Moser system. These new trigonometric R-matrices were
obtained from elliptic ones by applying to them some singular gauge transformation which
depends  on modular parameter $k$ and singular at $k=0$.
Thus the singular gauge transformation introduced in this section is its classical
version.
\section{Conclusions}
Let us sum up shortly main results of the paper. We have demonstrated that the sl(2,$\mathbb{C}$)-tops are equivalent
 to the two-particle
 systems of Calogero-Moser type. We proved that a Hamiltonian of generic
top is equivalent to one of the three canonical Hamiltonians. We denoted these Hamiltonians as
$XYZ$, $XXZ^{'}$, and $XXX^{'}$ respectively. For each canonical Hamiltonian
we found natural Darboux coordinates on coadjoint orbits (the bosonisation formulas). In terms of these coordinates
dynamics of $XYZ$, $XXZ^{'}$, and $XXX^{'}$ tops is presented as dynamics of elliptical, trigonometrical or rational
 two-particle Calogero-Moser system respectively. In this way we showed the equivalence between the sl(2,$\mathbb{C}$)-tops
 and the two-particle  Calogero-Moser systems explicitly.\\
\\
\\
\begin{Large}\textbf{Acknowledgments} \end{Large}\\
\\
The author is grateful to M.Olshanetsky, A.Levin, and A.Zotov for fruitful discussions and interest
to this work. The work was partly supported by RFBR grant 06-02-17382,   RFBR grant 06-01-92054-$KE_{a}$
and grant for support of scientific schools NSh-8004.2006.2.\\
\\

\end{document}